\documentclass[reqno]{amsart}
\usepackage{amssymb}
\usepackage[hypertex]{hyperref}
\allowdisplaybreaks[4]
\theoremstyle{plain}
\newtheorem{thm}{Theorem}
\newtheorem{lem}{Lemma}
\theoremstyle{remark}
\newtheorem{rem}{Remark}
\DeclareMathOperator{\td}{d\mspace{-1mu}}
\date{Drafted on 7 April 2008 in Melbourne}
\date{}

\begin{document}

\title{An alternative proof of Elezovi\'c-Giordano-Pe\v{c}ari\'c's theorem}

\author[F. Qi]{Feng Qi}
\address[F. Qi]{Research Institute of Mathematical Inequality Theory, Henan Polytechnic University, Jiaozuo City, Henan Province, 454010, China}
\email{\href{mailto: F. Qi <qifeng618@gmail.com>}{qifeng618@gmail.com}, \href{mailto: F. Qi <qifeng618@hotmail.com>}{qifeng618@hotmail.com}, \href{mailto: F. Qi <qifeng618@qq.com>}{qifeng618@qq.com}}
\urladdr{\url{http://qifeng618.spaces.live.com}}

\author[B.-N. Guo]{Bai-Ni Guo}
\address[B.-N. Guo]{School of Mathematics and Informatics, Henan Polytechnic University, Jiaozuo City, Henan Province, 454010, China}
\email{\href{mailto: B.-N. Guo <bai.ni.guo@gmail.com>}{bai.ni.guo@gmail.com}, \href{mailto: B.-N. Guo <bai.ni.guo@hotmail.com>}{bai.ni.guo@hotmail.com}}
\urladdr{\url{http://guobaini.spaces.live.com}}

\begin{abstract}
In the present note, an alternative proof is supplied for Theorem~1 in [N.~Elezovi\'c, C.~Giordano and J.~Pe\v{c}ari\'c, \textit{The best bounds in Gautschi's inequality}, Math. Inequal. Appl. \textbf{3} (2000), 239\nobreakdash--252.].
\end{abstract}



\thanks{The first author was partially supported by the China Scholarship Council}

\thanks{This paper was typeset using \AmS-\LaTeX}

\maketitle

\section{Introduction}

Let $s$ and $t$ be real numbers with $t-s\ne\pm1$. For $x\in(-\alpha,\infty)$, define
\begin{align}\label{psidef}
z_{s,t}(x)=\begin{cases}
\bigg[\dfrac{\Gamma(x+t)}{\Gamma(x+s)}\bigg]^{1/(t-s)}-x,&s\ne t,\\
e^{\psi(x+s)}-x,&s=t,
    \end{cases}
\end{align}
where $\alpha=\min\{s,t\}$,
\begin{equation}
\Gamma(x)=\int_0^\infty e^{-t}t^{x-1}\td t
\end{equation}
for $x>0$ stands for the classical Euler's gamma function, and $\psi(x)$ denotes the psi or digamma function, the derivative of the logarithm $\ln\Gamma(x)$.
\par
In order to bound the ratio of two gamma functions from both sides, N.~Elezovi\'c, C.~Giordano and J.~Pe\v{c}ari\'c proved in~\cite[Theorem~1]{egp} the following monotonicity and convexity results of the function $z_{s,t}(x)$.

\begin{thm}\label{egpthm}
The function $z_{s,t}(x)$ is either convex and decreasing for $\vert t-s\vert<1$ or
concave and increasing for $\vert t-s\vert>1$.
\end{thm}

The explicit or implicit origins and background of this theorem may be traced back to~\cite{gaut, kershaw, waston, wendel} and~\cite[Theorem~2]{Lazarevic}. This theorem or its special cases have been proved several times by different approaches in, for example, \cite{201-05-JIPAM, Lazarevic, notes-best-simple-open.tex, notes-best-simple.tex, notes-best.tex, notes-best.tex-rgmia, waston}. For detailed information on its history, please refer to the survey article~\cite{bounds-two-gammas.tex} published as a preprint recently.
\par
The purpose of this note is to supply an alternative proof for Theorem~\ref{egpthm}.

\section{Lemmas}

In order to prove Theorem~\ref{egpthm} alternatively, the following lemmas are necessary.

\begin{lem}[{\cite[p.~16]{magnus}}]\label{pgint}
The polygamma functions $\psi^{(n)}(x)$ can be expressed for $x>0$ and $n\in\mathbb{N}$ as
\begin{equation}
\psi ^{(n)}(x)=(-1)^{n+1}\int_{0}^{\infty}\frac{t^{n}}{1-e^{-t}}e^{-xt}\td t. \label{psin}
\end{equation}
\end{lem}

\begin{lem}[\cite{LaplaceTransform.html}]\label{convlotion}
Let $f_i(t)$ for $i=1,2$ be piecewise continuous in arbitrary finite
intervals included on $(0,\infty)$, suppose there exist some constants
$M_i>0$ and $c_i\ge0$ such that $\vert f_i(t)\vert\le M_ie^{c_it}$ for $i=1,2$.
Then
\begin{equation}
\int_0^\infty\bigg[\int_0^tf_1(u)f_2(t-u)\td u\bigg]e^{-st}\td t
=\int_0^\infty f_1(u)e^{-su}\td u\int_0^\infty f_2(v)e^{-sv}\td v.
\end{equation}
\end{lem}

\begin{lem}\label{qstcov}
For $u\in\mathbb{R}$ and $\beta>\alpha\ge0$ with $(\alpha, \beta)\ne(0,1)$, let
\begin{equation}\label{qfdef}
q_{\alpha,\beta}(u)=\begin{cases}
\dfrac{e^{-\alpha u}-e^{-\beta u}}{1-e^{-u}},&u\ne0;\\
\beta-\alpha,&u=0.
\end{cases}
\end{equation}
\begin{enumerate}
\item
The function $q_{\alpha,\beta}(u)$ is logarithmically convex for $\beta-\alpha>1$ and logarithmically concave for $0<\beta-\alpha<1$ on $(-\infty,\infty)$.
\item
For $\beta-\alpha>1$, the function
\begin{equation}\label{Qs,t;lambda(u)}
Q_{s,t;\lambda}(u)=q_{\alpha,\beta}(u)q_{\alpha,\beta}(\lambda-u)
\end{equation}
is increasing on $\bigl(\frac\lambda2, \infty\bigr)$ and decreasing on $\bigl(-\infty, \frac\lambda2\bigr)$, where $\lambda$ is any real constant; For $0<\beta-\alpha<1$, it is decreasing on $\bigl(\frac\lambda2, \infty\bigr)$ and increasing on $\bigl(-\infty, \frac\lambda2\bigr)$.
\end{enumerate}
\end{lem}

\begin{proof}
It is clear that the function $q_{\alpha,\beta}(u)$ can be rewritten as
\begin{equation*}
q_{\alpha,\beta}(u)=\frac{\sinh((\beta-\alpha)u/2)}{\sinh(u/2)} \exp\frac{(1-\alpha-\beta)u}2\triangleq p_{\alpha,\beta}\biggl(\frac{u}2\biggr).
\end{equation*}
Since the functions $q_{\alpha,\beta}(u)$ and $p_{\alpha,\beta}(u)$ are positive for $\beta>\alpha$, taking the logarithm of $p_{\alpha,\beta}(u)$ and differentiating yield
\begin{align*}
\ln p_{\alpha,\beta}(u)&=\ln\sinh((\beta-\alpha)u)-\ln\sinh u+(1-\alpha-\beta)u,\\
[\ln p_{\alpha,\beta}(u)]'&=(\beta-\alpha)\coth((\beta-\alpha)u) -\coth u-\alpha-\beta+1,\\
[\ln p_{\alpha,\beta}(u)]''&=\frac1{u^2}\biggl\{\biggl(\frac{u}{\sinh u}\biggr)^{2} -\biggl[\frac{(\beta-\alpha)u}{\sinh((\beta-\alpha)u)}\biggr]^{2}\biggr\}\\ &\triangleq\frac{[h(u)]^2-[h((\beta-\alpha)u)]^2}{u^2}.
\end{align*}
It is clear that the functions $h(u)$ and $[\ln p_{\alpha,\beta}(u)]''$ are even and the former is positive on $(-\infty,\infty)$, increasing on $(-\infty,0)$ and decreasing on $(0,\infty)$. As a result,
\begin{enumerate}
\item
for $\beta-\alpha>1$, if $u>0$, then $(\beta-\alpha)u>u>0$ and $h((\beta-\alpha)u)<h(u)$, and so $[\ln p_{\alpha,\beta}(u)]''>0$ on $(0,\infty)$;
\item
for $\beta-\alpha>1$, if $u<0$, then $(\beta-\alpha)u<u<0$ and $h((\beta-\alpha)u)<h(u)$, and so $[\ln p_{\alpha,\beta}(u)]''>0$ on $(-\infty,0)$;
\item
for $0<\beta-\alpha<1$, if $u>0$, then $0<(\beta-\alpha)u<u$ and $h((\beta-\alpha)u)>h(u)$, and so $[\ln p_{\alpha,\beta}(u)]''<0$ on $(0,\infty)$;
\item
for $0<\beta-\alpha<1$, if $u<0$, then $0>(\beta-\alpha)u>u$ and $h((\beta-\alpha)u)>h(u)$, and so $[\ln p_{\alpha,\beta}(u)]''<0$ on $(-\infty,0)$.
\end{enumerate}
From the obvious relationship $p_{\alpha,\beta}(u)=q_{\alpha,\beta}(2u)$ on $(-\infty,\infty)$, the logarithmically convex properties in Lemma~\ref{qstcov} follows readily.
\par
Taking the logarithm of $Q_{s,t;\lambda}(u)$ and differentiating give
\begin{equation*}
[\ln Q_{s,t;\lambda}(u)]'=\frac{q_{\alpha,\beta}'(u)}{q_{\alpha,\beta}(u)} -\frac{q_{\alpha,\beta}'(\lambda-u)}{q_{\alpha,\beta}(\lambda-u)}.
\end{equation*}
For $\beta-\alpha>1$, by the logarithmic convexities of $q_{\alpha,\beta}(u)$, it follows that the function $\frac{q_{\alpha,\beta}'(u)}{q_{\alpha,\beta}(u)}$ is increasing and $\frac{q_{\alpha,\beta}'(\lambda-u)}{q_{\alpha,\beta}(\lambda-u)}$ is decreasing on $(-\infty,\infty)$; From the obvious fact that $[\ln Q_{s,t;\lambda}(u)]'|_{u=\lambda/2}=0$, it follows that $[\ln Q_{s,t;\lambda}(u)]'>0$ for $u>\frac{\lambda}2$ and $[\ln Q_{s,t;\lambda}(u)]'<0$ for $u<\frac{\lambda}2$; Hence, the function $Q_{s,t;\lambda}(u)$ is increasing for $u>\frac\lambda2$ and decreasing for $u<\frac\lambda2$. Similarly, for $0<\beta-\alpha<1$, the function $Q_{s,t;\lambda}(u)$ is decreasing for $u>\frac\lambda2$ and increasing for $u<\frac\lambda2$. The proof of Lemma~\ref{qstcov} is proved.
\end{proof}

\begin{lem}\label{comp-thm-1}
For $x\in(0,\infty)$,
\begin{equation}\label{psi'ineq}
\ln x-\frac1x<\psi(x) <\ln x-\frac1{2x}
\end{equation}
and
\begin{equation}\label{psiineq}
\frac1{x}+\frac1{2x^2}< \psi'(x)<\frac1{x}+\frac1{x^2}.
\end{equation}
\end{lem}

\begin{proof}
This may be derived easily from the fact~\cite[p.~82]{e-gam-rat-comp-mon} that a completely monotonic function which is non-identically zero cannot vanish at any point on $(0,\infty)$ and the complete monotonicity obtained in~\cite[Theorem~2]{sandor-gamma-2-ITSF.tex} and~\cite[Theorem~2]{sandor-gamma-2-ITSF.tex-rgmia}: The function $\psi(x)-\ln x+\frac{\alpha}x$ is completely monotonic on $(0,\infty)$ if and only if $\alpha\ge1$ and so is the function $\ln x-\frac{\alpha}x-\psi(x)$ if and only if $\alpha\le\frac12$.
\end{proof}

\section{An alternative proof of Theorem \ref{egpthm}}

Since $z_{s,t}(x)=z_{t,s}(x)$, without loss of generality, we can assume $t>s\ge0$ and $t-s\ne1$ in what follows.
\par
Differentiation of $z_{s,t}(x)$, utilization of~\eqref{psin} and application of Lemma~\ref{convlotion} yield
\begin{gather}
z'_{s,t}(x)=\frac{[z_{s,t}(x)+x][\psi(x+t)-\psi(x+s)]}{t-s}-1,\label{zi}\\
\frac{z''_{s,t}(x)}{z_{s,t}(x)+x}=\biggl[\frac{\psi(x+t)-\psi(x+s)}{t-s}\biggr]^2 +\frac{\psi'(x+t)-\psi'(x+s)}{t-s}\notag\\
=\biggl[\frac1{t-s}\int_s^t\psi'(x+u)\td u\biggr]^2+\frac1{t-s}\int_s^t\psi''(x+u)\td u\notag\\
=\biggl[\frac1{t-s}\int_s^t\int_{0}^{\infty}\frac{ve^{-(x+u)v}}{1-e^{-v}}\td v\td u\biggr]^2 -\frac1{t-s}\int_s^t\int_{0}^{\infty}\frac{v^2e^{-(x+u)v}}{1-e^{-v}}\td v\td u\notag\\
=\biggl(\int_{0}^{\infty}\frac{ve^{-xv}}{1-e^{-v}}\cdot\frac1{t-s}\int_s^te^{-uv}\td u\td v\biggr)^2 -\int_{0}^{\infty}\frac{v^2e^{-xv}}{1-e^{-v}}\cdot\frac1{t-s}\int_s^te^{-uv}\td u\td v\notag\\
=\biggl(\int_{0}^{\infty}\frac{e^{-xv}}{1-e^{-v}}\cdot\frac{e^{-s v}-e^{-t v}}{t-s}\td v\biggr)^2 -\int_{0}^{\infty}\frac{ve^{-xv}}{1-e^{-v}}\cdot\frac{e^{-s v}-e^{-t v}}{t-s}\td v\notag\\
=\int_{0}^{\infty}\biggl[\frac1{(t-s)u} \int_{0}^{u}q_{s,t}(r)q_{s,t}(u-r)\td r -q_{s,t}(u)\biggr]ue^{-xu}\td u\notag\\ \label{line-last}
=\int_{0}^{\infty}\biggl[\frac1{(t-s)u} \int_{0}^{u}Q_{s,t;u}(r)\td r-q_{s,t}(u)\biggr]ue^{-xu}\td u.
\end{gather}
If $t-s>1$, by the monotonicity of $Q_{s,t;\lambda}(u)$ in Lemma~\ref{qstcov}, it follows easily that
$$
Q_{s,t;u}(r)\le Q_{s,t;u}(0)=Q_{s,t;u}(u)=q_{s,t}(0)q_{s,t}(u)=(t-s)q_{s,t}(u),
$$
consequently, the bracketed term in the line~\eqref{line-last} is negative on $(0,\infty)$, and so $z''_{s,t}(x)<0$. If $0<t-s<1$, the similar argument leads to $z''_{s,t}(x)>0$. The convex and concave properties of $z_{s,t}(x)$ are proved.
\par
By the mean value theorem, it is immediate that
\begin{align*}
z'_{s,t}(x)+1&=\bigg[\bigg(\frac{\Gamma(x+t)}{\Gamma(x+s)}\bigg)^{1/(t-s)} \frac{\psi(x+t)-\psi(x+s)}{t-s}\bigg]\\
&=\frac{\psi(x+t)-\psi(x+s)}{t-s} \exp\frac{\ln\Gamma(x+t)-\ln\Gamma(x+s)}{t-s}\\
&=\psi'(x+\xi_1)e^{\psi(x+\xi_2)},\quad \text{$\xi_i\in(s,t)$ for $i=1,2$}.
\end{align*}
By inequalities in~\eqref{psi'ineq} and~\eqref{psiineq}, it is ready to obtain
\begin{equation*}
\biggl[\frac{x+\xi_2}{x+\xi_1}+\frac{x+\xi_2}{2(x+\xi_1)^2}\biggr]\frac1{e^{1/(x+\xi_2)}} <z'_{s,t}(x)+1 <\biggl[\frac{x+\xi_2}{x+\xi_1}+\frac{x+\xi_2}{(x+\xi_1)^2}\biggr]\frac1{e^{1/2(x+\xi_2)}}
\end{equation*}
which means $\lim_{x\to\infty}z'_{s,t}(x)=0$. For $t-s>1$, the conclusion that $z''_{s,t}(x)\le0$ obtained above implies $z'_{s,t}(x)$ is decreasing, and so $z'_{s,t}(x)>0$ and $z_{s,t}(x)$ is increasing. For $0<t-s<1$, the result that $z''_{s,t}(x)\ge0$ obtained above implies $z'_{s,t}(x)$ is increasing, and so $z'_{s,t}(x)<0$ and $z_{s,t}(x)$ is decreasing. The proof of Theorem~\ref{egpthm} is complete.

\section{Some remarks}

\begin{rem}
The logarithmically convex properties of $q_{\alpha,\beta}(u)$ on $(-\infty,0)$ in Lemma~\ref{qstcov} of this paper corrects some mistakes appeared in~\cite[Lemma~1]{notes-best.tex} and~\cite[Lemma~1]{notes-best.tex-rgmia}. However, these mistakes did not affect the correctness of the proof provided in~\cite{notes-best.tex, notes-best.tex-rgmia} for Theorem~\ref{egpthm}, since properties of $q_{\alpha,\beta}(u)$ on $(-\infty,0)$ are idle there.
\end{rem}

\begin{rem}
The logarithmically convex properties in Lemma~\ref{qstcov} of this paper were also proved in~\cite{mon-element-exp.tex-rgmia} by using different techniques. Also see~\cite{mon-element-exp-AIMS.tex} and related references therein.
\end{rem}

\begin{rem}
It is well-known that a positive and $k$-times differentiable function $f(x)$ is said to be $k$-log-convex (or $k$-log-concave, respectively) on an interval $I$ with $k\ge2$ if and only if $[\ln f(x)]^{(k)}$ exists and $[\ln f(x)]^{(k)}\ge0$ (or $[\ln f(x)]^{(k)}\le0$, respectively) on $I$. The $3$-log-convex properties of $q_{\alpha,\beta}(u)$ were already obtained in~\cite[Theorem~1.1]{comp-mon-element-exp.tex}: For $1>\beta-\alpha>0$, the function $q_{\alpha,\beta}(u)$ is $3$-log-convex on $(0,\infty)$ and $3$-log-concave on $(-\infty,0)$; For $\beta-\alpha>1$, it is $3$-log-concave on $(0,\infty)$ and $3$-log-convex on $(-\infty,0)$.
\end{rem}

\end{document}